\documentclass[12pt,a4paper]{article} 
\usepackage[utf8]{inputenc} 
\usepackage[english]{babel} 
\usepackage{amsmath} 
\usepackage{amsfonts} 
\usepackage{amssymb} 
\usepackage{makeidx}
\usepackage{graphicx}
\usepackage[left=2cm,right=2cm,top=2cm,bottom=2cm]{geometry}
\author{Ramachandra BHAT \\Rajajinagar, Bengaluru, \\Karnataka, India - 560010.}
\title{Locus of Intersection for Trisection}

\begin{document} 
\maketitle 

\begin{abstract}
While solving problems, if direct methods does not provide solution, indirect methods are explored.  Today, we need an indirect method to solve the problem of angle trisection as the direct methods have been proved using algebra to be impossible by Euclidean Geometry, using only straight edge and compasses.   The unstoppable curiosity of Geometers and the newer advanced tools available with time have led to newer approaches to progress further.  There is a method using Origami (paper folding) procedure to trisect an angle.  The algebraic analysis of this procedure gives us a method of finding trisection using a locus of a point of intersection of two circles.  The algebraic  analysis and the equation for the locus of the point of intersection of two circles leading to trisection of any given unknown angle without any measurements is described here. 

\end{abstract}

\section{Introduction} 

In ancient times, there were no formal measurement systems in place.  While working on designs, Greeks defined the dimensions by considering an arbitrary length as a unit of length, for that particular design.  They derived other dimensions within that design geometrically by the addition, multiplication, subtraction and division operations.  In addition, they knew the right angle and the Pythagoras theorem that helped them to get the square roots of the given line segment.  Therefore, in Euclidean Geometry, one does not measure the dimensions but the expected designs were very precise, accurate and simple to reproduce.  The only geometrical tools used in those days were a ruler (unmarked straightedge) and Compasses.  

With the advancement of mathematics and allied subjects, today we could actually make the measurements to show the accuracy and precision of the Greeks and explain them through the algebraic and trigonometric formulations.  However, even today, the Three Famous Problems of Geometry have remained unsolved.  One of them being the trisecting an angle by the Euclidean procedure.   

While many Geometers have been exploring to find solutions to these problems, in 1837, Pierre Wantzel proved that finding a solution to the problem of trisection of an angle is impossible\cite{Wiki, Jagadeeshan}.  This only meant, using only ruler and compasses, it is not possible to trisect an angle of any given value, with the \textit{tools and knowledge available at that instant of time}. It should also be noted that the proof of impossibility considers primarily the constructability of the angle of value equivalent to one-third of the given value and not the trisectability of the given angle of unknown measure directly. 

In the literature,  there are various methods described \cite{Yates, Mallik, Brooks, Forum, Wolfram} to trisect an angle using various additional tools (Marked ruler, trisector, Quadratrix of Hippias, hyperbola, tomahawk, linkages, etc.) and hence, are not Euclidean methods.  However, the only exact and simple procedure for the trisection of any given angle of unknown measure demonstrated as well as proved algebraically is by Origami \cite{Richeson, Shima, Kung}, i.e., the Japanese paper folding technique, by Hisashi Abe.  However, in this procedure, multiple conditions get satisfied simultaneously in one folding operation, by using an implicit marked ruler.  The results of the algebraic analysis of this Origami method are utilized in exploring the trisection using ruler and compasses only.  While doing so, it has been observed that the locus of a point of intersection of two circle gives the exact results for trisection, though it would take a few extra steps to generate the locus graph to implement the procedure by Euclidean Geometry. 

\section{The problem definition:} 
The two approaches followed by many geometers in their exploration are: \begin{enumerate}\item[] (i) to divide the given angle into three equal parts, and/or \item[] (ii) to explore the constructability of an angle of value equivalent to one-third of the given angle. 
\end{enumerate}  
 
For example, some select angles such as 0, 45, 72, 90, 108 and 180 degrees were trisectable \cite{Yates2}, by directly constructing a corresponding angle equivalent to their one-third value and not by trisecting the given angle. This means, \begin{enumerate} \item[(i)] the actual value of the given angle was known (or measured) prior to attempting the trisection, and \item[(ii)] its trisected component is constructible by Euclidean geometry.  \end{enumerate}  
 
Since the measure of the given angle was unknown, Greeks could only explore the trisectability and not constructibility.  In the present study, the intention was to trisect the given angle and not to find the constructability of angle of any given value.  Hence, the definition by the Euclidean geometry may be presented as: 

\emph{\textbf{``Divide the given angle} (of unknown value) \textbf{into three equal angles} (angular parts) using only two tools, viz., (i) an unmarked straight edge (ruler) and (ii) a compasses".} 

\section{Understanding Trisection} 

Let us consider the addition rule for Sin(A+B).  Let A = $\theta$ and B = $2*\theta$ and draw the diagram for the same.  Next, consider A = $2*\theta$ and B = $\theta$ and draw the corresponding diagram with the same unit length and angle value $\theta$.  Overlay the two diagrams (on same scale).   This is presented in an overladen square ABCD in Fig. 1.  
 
\begin{figure}[ht!]
\includegraphics[scale=0.5]{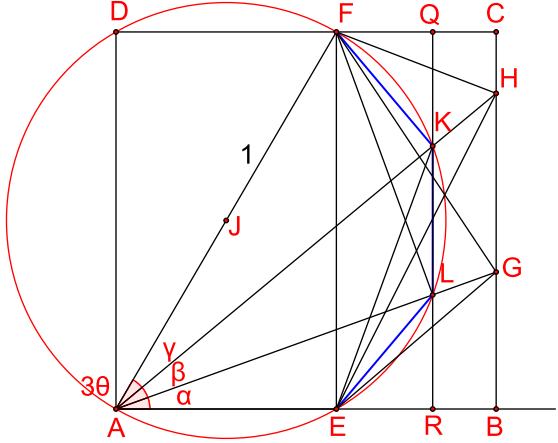} 
\caption{Trigonometric addition formula diagram and Trisection} 
\label{Trisection Analysis using Addition Formula} 
\end{figure}  
 
Let the given angle $\angle$BAF = $3\theta$. 
When this angle is trisected, let us consider the three equal components as $\alpha, \beta \ and \ \gamma$, i.e., $3\theta = \alpha + \beta + \gamma$. 

In addition, let us draw a circle with radius = 0.5 unit of length with J (= mid point of AF) as the center.  K and L are the intersection points of the circle with AH and AG, the lines of trisected angles.  Now, from the diagram, it could be found that segments FK = KL = LE = $Sin(\theta)$. 

Hence, as exprected, the angle trisection corresponds to the trisection of an arc of a circle and the length of cords corresponding to each section is given by sine of the trisected component of the given angle [$sin(\theta$)].     
 
\begin{center} 
\begin{tabular}{c}  
From the diagram, it may be noted that \\
\hline 
$\angle BAF = 3\theta$ and $3\theta = \alpha + \beta + \gamma$ and  $\alpha = \beta = \gamma$ \\ 
$\angle BAG = \alpha$; $\angle GAH = \beta$; $\angle HAF = \gamma$. \\ 
AE = $Cos(3\theta)$; EF = $Sin(3\theta)$; AF = 1 = AG. \\ 
\hline  
BG = $Sin(\theta)$; GF = $2*Sin(\theta)$; AF = 1 = AG. \\ 
In the Origami procedure, BG corresponds to the first fold. \\
\hline 
\end{tabular} 
\end{center} 
 
\section{Understanding Origami method} 

The geometrical picture of the Origami procedure is given in Fig. 2.

\begin{figure}[ht!]
\includegraphics[scale=0.4]{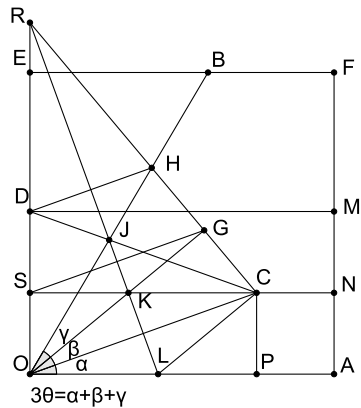} 
\caption{Hisashi Abe's Origami diagram for Trisection}  
\label{Origami Procedure diagram}  
\end{figure} 

Let us assume the given angle $\angle$AOB = 3$\theta$ (for convenience).  Considering the similar isosceles triangles DCO and HOC, it has been proved already that the three angles $\alpha, \beta$ and $\gamma$ are all mutually equal \cite{Richeson}.  Hence, $3\theta = \alpha + \beta + \gamma = 3*\theta$, i.e., to say, $\alpha = \beta = \gamma = \theta$, and hence the trisection.  
 
Here, segments OH = OC and OS = SD = HG = GC = CP.  

If we define OH = OC as the unit of length (=1),  then OS = CP = sin($\theta$) and OP = cos($\theta$).  

From the Origami method, it is clear that the two unknowns in the procedure for trisection are the unit of length (= 1) and the fold length (= sin($\theta$)) and they are interdependent.  While folding the paper, the value for the fold length (= OS = CP = sin($\theta$)) is arbitrarily assumed and the corresponding value for the unit length (= 1 = OH) is determined. While doing so, the segment OD is matched with CH such that OH = OC = 1.  Here, the segment OD serves as the \emph{marked ruler} to identify the location of point G (corresponding to the point S on OD) which is the advantage over the Euclidean method.  
 
\section{Options for solution}  

From the Origami procedure and the addition formula for the trigonometric function [Sin(A+B)], it is clear that we need to explore a method to find the value of either unit length (segment length=1) or the value for sin($\theta$) to achieve the trisection, while we could assume the value of one of these two interdependent parameters arbitrarily.  Now, the options available for the exploration of trisection of a given angle is summarized in the following GFR (Given-Find-Report) table.  

\begin{center}  
\begin{tabular}{||l|c|c|c||} 
\hline 
\multicolumn{4}{||c||}{GFR Table for the solution exploration} \\
\hline  
Options & Given \{\& Assumed\} & Parameter to Find & Report \\
\hline 
First  & $3\theta$, Paper, \{$\& Sin\theta$\} & 1 (or Cos$\theta$) & Origami procedure \\
\hline 
Second & 3$\theta$, S \& C, \{\& Sin$\theta$\} & 1 (or Cos$\theta$) & 	Euclidean procedure* \\ 
\hline 
Third  & 3$\theta$, S \& C, \{\& 1\} & Sin$\theta$ (or Cos$\theta$) &	Euclidean procedure* \\
\hline 
Fourth & 3$\theta$, S \& C, \{\& Sin$\theta$\} & 1 (or $Cos\theta$) & Locus method** \\
\hline 
\multicolumn{4}{||c||}{* Using only S \& C (Straightedge and Compassess); ** Present study} \\
\hline  
\end{tabular} 
\end{center}  

\section{Locus of the point of intersection} 

In the present study, the locus of a point of intersection of two circles is made use of to find the exact dimension of the unit length segment, while assuming a value for $sin(\theta)$.   

\subsection{\textbf{Construction}}

Let the given angle be $\angle$AOB = 3$\theta$ (= 3*$\theta$). 
Similar to the Origami procedure, let us mark the two folds as lines parallel to OA (base segment of the given angle) at C and D.  The length OC = CD corresponds to Sin$\theta$ (assumed value) for the given angle.  

\begin{figure}[ht!]
\includegraphics[scale=0.33]{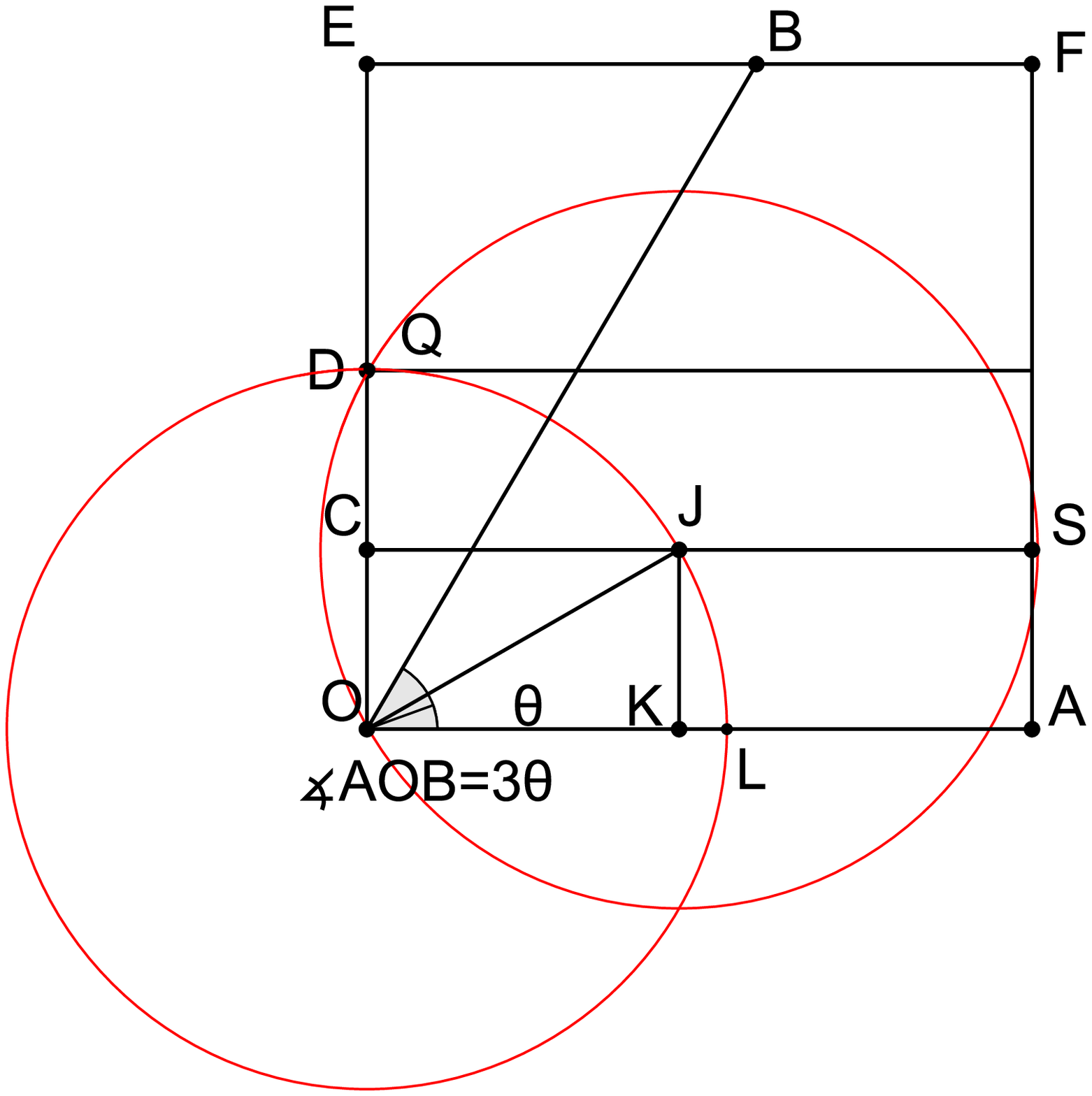} 
\caption{Diagram for the locus of the intersection}  
\label{Locus Definition diagram} 
\end{figure} 

Now, to find the corresponding value of the unit length, draw the first circle with O as the center and OD as the radius.  Let the circle cut the first fold line at J.  Mark the point of projection of point J on OA as K.  At this stage, KJ = OC and OJ = OD = 2 * OC.  Hence, $\angle$KOJ = 30$^\circ$ and OK = $\sqrt{3}$*KJ. 

At J, draw the second circle with J as the center and 2*JK (= OD) as the radius.  The point of intersection of this second circle with the first circle is Q.  To begin with, the point Q coincides with D.  
 
Now, to get the locus of the point Q, vary the position of J along the line CS such that the length OJ increases gradually, till the locus line crosses the line OB.  Let the point of intersection of the locus line with the line OB is N.  Then, ON = 1 (unit of length).  When ON = 1, the length OJ is also of unit length (radius of the same circle) and JK = sin($\theta$), where $\theta$ is the value of the trisected angle, i.e., $\angle$AOB = 3*$\angle$AOJ.  Here, JN = 2 * KJ = OD, same as in Origami procedure.  Hence, the proof of trisection is similar to that for the Origami procedure.  

\begin{figure}[ht!] 
\includegraphics[scale=0.30]{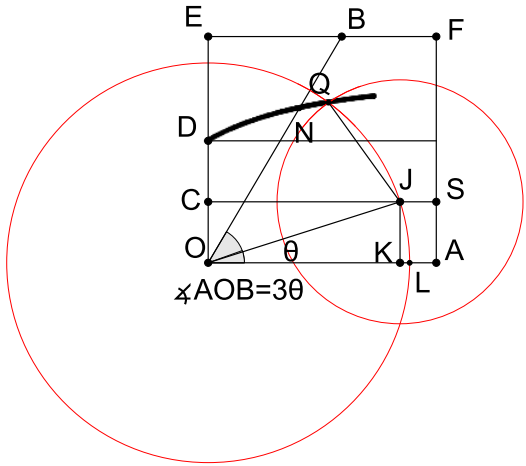} 
\caption{Locus of intersection and the trisection}  
\label{Locus of Intersection for trisection} 
\end{figure} 

\subsection{Equation for the Locus of Point Q} 

Let OC = JK = a and OK = b.  The coordinates of the points are O = (0, 0); D = (0, 2a); J = (b, a).

$\therefore$ The equation for Circle 1: $(x^2 + y^2) = OJ^2 = OQ^2 = (b^2 + a^2)$. 
  
and the equation for Circle 2:  $(x-b)^2 + (y-a)^2) = JQ^2 = (2a)^2 = 4a^2)$. 

At Q - the point of intersection, 

 $(x^2 + b^2 - 2bx) + (y^2 + a^2 - 2ay) = 4a^2$. 
 
 $(x^2 + y^2) - 2(bx + ay) + (b^2 + a^2) = 4a^2$. 

 $(b^2 + a^2) - 2(bx + ay) + (b^2 + a^2) = 4a^2$. 
 
 $\therefore$ $(bx + ay) = (b^2 - a^2)$  or 

 $ y = (\frac{-b}{a})*x + \frac{(b^2 - a^2)}{a}$. 
 
The value of b is variable (not constant) and varies with the values of x and y.  Hence, the equation for the locus corresponds to a curve and not a straight line.  Every point Q on the locus curve corresponds to $\angle$AOQ (= $3\delta$) for which $\angle$AOJ = $\delta$ is its trisected angle.  Hence,  the procedure gives the trisection for the spectrum of angles from $\angle$AOD (= 90$^\circ$) till $\angle$AOB (= given angle = $3\theta$).  Therefore, when Q lies on the line OB, it gives the trisection for the given angle $\angle$AOB.  

When Q coincides with D (in the beginning), x = 0, y = 2a, b = $\sqrt{3}$*a.  

When Q coinsides with N (at the trisection), 

a = $sin(\theta$) and b = $cos(\theta$).  Then, x = $cos(3\theta$) and y = $sin(3\theta$).  

Hence, the value of b varies between $\sqrt{3}$*a and $cos(\theta)$ while the value of x varies between 0 and $cos(3\theta)$ for the locus of the point Q between point D and Point N.   

\section{Summary and Conclusions} 
 
The proof for the trisection of a given unknown angle by the Origami method (Abe's method) is well established.  Using the same concept of the dimensions and the intersections, trisection procedure using straight edge and compass is explored in the present study.  The locus of the point of intersection of two circles establishes the trisection procedure.  To draw the complete locus line, we need to find the points at various locations using compass.  Once the concept is understood, it would be easy to do so.  When the locus line is very close to the intersection point with the angle line (OB), it is expected to be linear at that very short length.  Hence, by finding two points, one before and one after the line (OB), it is possible to find the exact unit length and hence the trisection.  This locus graph is a universal plot and can be made separately and used for any angle between 0 and 90$^\circ$.  Hence, it is possible, not impossible, to trisect an angle with only unmarked ruler and compasses.  It only needs extra effort to find the exact point of intersection of the line of locus with the the angle line (OB). 

\section{Acknowledgment}  

The support and cooperation of my family during the exploration is gratefully appreciated.  
In the early stage of this study (Mar-2016), while exploring suitable software to verify the accuracy of the ideas, the suggestion of Mr. John Page (mathopenref@gmail.com) to use GeoGebra \cite{Geogebra} was of great help and the same is acknowledged here. 

All the explorations of different options were done using the GeoGebra software.  Hence, a special thanks for the creator of the GeoGebra software who made it so user-friendly and available freely on the Internet \cite{Geogebra}.

\subsection*{Conflict of Interest}
There is no conflict of interest.


\begin{center}  
\begin{tabular}{lcl}
\hline 
Dr Ramachandra Bhat,  & \hspace{50pt} & \textit{Alumni of}\\
115, 12$^{th}$ B Main Road, &  & 1. R V High School, Itgi, \\ 
6th Block, Rajajinagar, &  &   Taluk: Siddapur (Uttara Kannada), \\
Bengaluru - 560 010. &  & Karnataka, INDIA.  \\
Karnataka, INDIA. &  & 2. Sarada Vilas College, Mysore.\\
Mobile: +91-9902461175. &  & 3. Yuvaraja's College, Mysore. \\ 
Email: $rs\_bhat@yahoo.com $ &  & 4. Indian Institute of Technology, Bombay.\\
\hline 
\end{tabular} 
\vspace{24pt}
\fbox{\textbf{Mathematicians are BORN to SOLVE Problems}} 
\end{center} 


\begin{thebibliography}{1}

\bibitem{Wiki} 
\textit{Angle Trisection}, https://en.wikipedia.org/wiki/Angle-trisection.  

\bibitem{Jagadeeshan} 
Shashidhar Jagadeeshan, \textit{Whoever said nothing is impossible? Three problems from Antiquity}, Resonance, (March 1999), 25--38.

\bibitem{Yates} 
Robert C Yates, \textit{Classics in Mathematics Education, Vol.3, The Trisection Problem}, Chapter I, The problem, 1942 (Reprinted 1971, National Council of Teachers of Mathematics). 

\bibitem{Mallik} 
AK Mallik, \textit{Kempe’s angle trisector: Sir Alfred Bray Kempe – An Amateur Kinematician}, Resonance, (Mar 2011), 204--218. (Angle Trisector, 214--215). 

\bibitem{Brooks} 
David Alan Brooks, \textit{A new method of trisection}, The college Mathematics Journal, \textbf{38(2)}, (March 2007), Mathematical Association of America, 78--81. 

\bibitem{Forum} 
\textit{Angle Trisection}, http://www.geom.uiuc.edu/docs/forum/angtri/

\bibitem{Wolfram} 
\textit{Angle Trisection}, http://mathworld.wolfram.com/AngleTrisection.html 

\bibitem{Richeson} 
Dave Richeson, \textit{Angle trisection using Origami }(Hisashi Abe’s trisection, 1970), posted on June 1, 2012, https://divisbyzero.com/2012/06/01/angle-trisection-using-origami/   

\bibitem{Shima} 
Hiroyuki Shima, \textit{Origami revisited}, Am. J Applied Math, (2013), 39--43. 

\bibitem{Kung} 
Sidney H Kung, \textit{An interesting example of angle trisection by paper folding}, \\
https://www.cut-the-knot.org/pythagoras/PaperFolding/KungAngleTrisection.shtml, Jan. 6, 2011. \\
\textit{Angle Trisection By Paper Folding}, \\https://www.cut-the-knot.org/pythagoras/PaperFolding/AngleTrisection.shtml

\bibitem{Yates2} 
Ref.\cite{Yates}, page 16. 

\bibitem{Geogebra} 
https://www.geogebra.org/download? lang=en. 
\end{thebibliography}
\end{document}